\documentclass[reqno]{amsart}

\newtheorem{theorem}{Theorem}
\newtheorem{lemma}[theorem]{Lemma}
\newtheorem{proposition}[theorem]{Proposition}
\newtheorem{corollary}[theorem]{Corollary}
\newtheorem{example}[theorem]{Example}
\theoremstyle{definition}

\theoremstyle{remark}

\numberwithin{equation}{section}

\newcommand{\intav}[1]{\mathchoice {\mathop{\vrule width 6pt height 3 pt depth  -2.5pt
\kern -8pt \intop}\nolimits_{\kern -6pt#1}} {\mathop{\vrule width
5pt height 3  pt depth -2.6pt \kern -6pt \intop}\nolimits_{#1}}
{\mathop{\vrule width 5pt height 3 pt depth -2.6pt \kern -6pt
\intop}\nolimits_{#1}} {\mathop{\vrule width 5pt height 3 pt depth
-2.6pt \kern -6pt \intop}\nolimits_{#1}}}

\newcommand{\intavl}[1]{\mathchoice {\mathop{\vrule width 6pt height 3 pt depth  -2.5pt
\kern -8pt \intop}\limits_{\kern -6pt#1}} {\mathop{\vrule width 5pt
height 3  pt depth -2.6pt \kern -6pt \intop}\nolimits_{#1}}
{\mathop{\vrule width 5pt height 3 pt depth -2.6pt \kern -6pt
\intop}\nolimits_{#1}} {\mathop{\vrule width 5pt height 3 pt depth
-2.6pt \kern -6pt \intop}\nolimits_{#1}}}

\def \N {\Bbb{N}}
\def \R {\Bbb{R}}

\begin{document}

\title[Maximal operators]{On the regularity of maximal operators}

\author[E. Carneiro]{Emanuel Carneiro}
\address{Department of Mathematics, University of Texas at Austin, Austin, TX 78712-1082.}
\email{ecarneiro@math.utexas.edu}
\thanks{The first author was supported by CAPES/FULBRIGHT grant BEX 1710-04-4.}

\author[D. Moreira]{Diego Moreira}
\address{Department of Mathematics, University of Iowa, Iowa City, IA 52242}
\email{dmoreira@math.uiowa.edu}

%    General info

\date{March, 11, 2008}


\begin{abstract}
We study the regularity of the bilinear maximal operator when applied to Sobolev functions, proving that it maps $W^{1,p}(\R) \times W^{1,q}(\R) \to W^{1,r}(\R)$ with $1 <p,q < \infty$  and $r\geq 1$, boundedly and continuously. The same result holds on $\R^n$ when $r>1$. We also investigate the almost everywhere and weak convergence under the action of the classical Hardy-Littlewood maximal operator, both in its global and local versions.
\end{abstract}

\maketitle

\section{Introduction}%%%%%%%%%%%%%%%%%%%%%%%%%%- Introduction%%%%%%%%%%%%%%%%%%%%%%%%%%%%%%%%%%%%%%%%%%%%%%%%%%%%%%%%%%

Over the last decade there has been considerable interest in understanding the regularity properties of maximal and singular integral operators, for instance how the weak differentiability is preserved. The first work in this direction is due to Kinnunen (\cite{K}) in 1997 when he observed that the classical Hardy-Littlewood maximal operator is bounded on the Sobolev space $W^{1,p}(\R^n)$ for $p>1$, using functional analytic tools (weak compactness arguments). Later developments on the subject include the boundedness of the local maximal operator for $W^{1,p}(\Omega),\, p>1$, where $\Omega$ is a proper open subset of $\R^n$ (\cite{KL}), and the continuity of the maximal operator for $W^{1,p}(\R^n),\, p >1$ (\cite{Lu}). Other interesting papers related to this topic are \cite{AP}, \cite{HO} and \cite{KS}.

In the first part of this paper we will be concerned with the following family of bilinear maximal operators in $\R^n$. For $\alpha \neq 1$ define
\begin{align}
\begin{split}
\mathcal{M}(f,g)(x) &= \sup_{R>0} \dfrac{1}{m(B_R)} \int_{B_R} |f(x-\alpha y)g(x-y)|dy  \label{BM1} \\
 &:=  \sup_{R>0}\intav{B_{R}} |f(x-\alpha y)g(x-y)|dy,
\end{split}
\end{align}
where $B_R$ is the ball of radius $R$ centered at the origin, and $m(A)$ denotes the $n$-dimensional Lebesgue measure of the measurable set $A\subset\mathbb{R}^{n}$. An application of H\"older's inequality tells us that this operator maps $L^p(\R^n) \times L^q(\R^n)$ into $L^r(\R^n)$ where $1/p + 1/q = 1/r$, $1<p,q < \infty$ and $r>1$. In 2000, M. Lacey in the remarkable paper \cite{L} showed that the family of one-dimensional bilinear maximal operators defined by (\ref{BM1}) maps $L^p(\R) \times L^q(\R)$ into $L^1(\R)$ where $1/p + 1/q = 1$, $1<p,q<\infty$, solving a conjecture posed by A. Calder\'on in 1964.

Throughout this paper we consider the following norm for $f\in W^{1,p}$:
$$ ||f||_{1,p} = ||f||_{p} + ||\nabla f||_{p} \  , $$ 
where $\nabla f$ denotes the weak gradient of the Sobolev function $f$. 

In the spirit of  Kinunnen's philosophy, M. Lacey's paper (\cite{L}) opens the way to raise the following natural question: what does happen when the bilinear maximal operator is applied to Sobolev functions? This brings our first result.

\begin{theorem}\label{thm1}
Given $\alpha \neq 1$, the bilinear maximal operator $\mathcal{M}$ defined in {\rm (\ref{BM1})} maps $W^{1,p}(\R^n) \times W^{1,q}(\R^n) \to W^{1,r}(\R^n)$ boundedly and continuously, where $1/p + 1/q = 1/r$,\ $1< p,q <\infty$ and
\begin{enumerate}
\item[(a)] $r \geq 1$, if $n=1$;
\item[(b)] $r>1$, if $n>1$.
\end{enumerate}
Boundedness is a consequence of the following pointwise estimate:
\begin{equation}\label{bound2}
|\nabla \mathcal{M}(f,g)(x)| \leq \mathcal{M}(f,|\nabla g|)(x) + \mathcal{M}(|\nabla f|,g)(x) \hspace{.3cm} \textrm{ a.e. x}  \in\mathbb{R}^{n}.
\end{equation}
\end{theorem}
Because of Lacey's theorem, the case $n=1$ and $r=1$ becomes the key difference between the bilinear maximal and the prior works on the classical Hardy-Littlewood maximal operator. The functional analytic arguments in \cite{K} and \cite{KL}, relying on the reflexivity of $L^r(\R^n)$ for $r>1$, are no longer available. To overcome this difficulty, we adopt here the approach introduced in \cite{HO}. For the continuity part we follow the insightful and elegant proof of Luiro in \cite{Lu}. Again the case $n=1$, $r=1$ is a new feature. We observe that the sublinearity of the maximal operator is not sufficient to infer continuity in $W^{1,p}$ from its boundedness.\\

{\it Remark}: It is believed that the bilinear maximal operator in $\R^n$, $n>1$, also maps $L^p(\R^n) \times L^q(\R^n)$ into $L^1(\R^n)$ if $1/p + 1/q = 1$,\ $1<p,q < \infty$. If this is indeed the case, we can include $r=1$, $n>1$ in Theorem \ref{thm1} above with our proof.\\

In a second and independent part of the paper (section 4), we study the behavior of the almost everywhere and weak convergence under the action of the classical Hardy-Littlewood maximal operator, both in its global and local versions.

\section{Proof of Theorem 1: Boundedness}%%%%%%%%%%%%%%%%%%%%%%%%%%%%%%-Section2%%%%%%%%%%%%%%%%%%%%%%%%%%%%%%%%%%%%%%%%%%

The proof of Theorem 1 is independent of the parameter $\alpha$ (as long as $\alpha \neq 1$ to guarantee Lacey's theorem) and henceforth we work with $\alpha = -1$. Let $f \in W^{1,p}(\R^n)$ and $g \in W^{1,q}(\R^n)$. Since $|f| \in W^{1,p}(\R^n)$ and $| \nabla |f| | = |\nabla f|$ we can assume that $f$ and $g$ are nonnegative. 

We start with $f,g \in C^{\infty}_0$ and fix $x,y \in \R^n$. We may assume that $\mathcal{M}(f,g)(x) \geq \mathcal{M}(f,g)(y)$. Let us take a sequence of radii $\{r_n\}_{n=1}^{\infty}$, $0< r_n < \infty$ such that
\begin{equation*}
 \lim_{n \to \infty} \intav{B_{r_n}} f(x+z)g(x-z)dz = \mathcal{M}(f,g)(x)
\end{equation*}
and write
\begin{equation*}
u_{r_n} (x) = \intav{B_{r_n}} f(x+z)g(x-z)dz 
\end{equation*}
for all $n \in \N$. Since
\begin{eqnarray*}
\left|\mathcal{M}(f,g)(x) - \mathcal{M}(f,g)(y)\right| &\leq& \left(\mathcal{M}(f,g)(x) - u_{r_n} (x)\right) + \left(u_{r_n} (x) - u_{r_n} (y)\right)
\end{eqnarray*}
for all $n \in \N$, we have
\begin{eqnarray}\label{limsup}
\left|\mathcal{M}(f,g)(x) - \mathcal{M}(f,g)(y)\right| &\leq& \limsup_{n\to \infty} \left(u_{r_n} (x) - u_{r_n} (y)\right).
\end{eqnarray}
By combining equation (\ref{limsup}) with the estimate
\begin{eqnarray*}
&& \left|u_{r_n} (x) - u_{r_n} (y)\right| \\
& = & \left|\intav{B_{r_n}} \Big\{ f(x+z)g(x-z)dz - f(y+z)g(y-z)\Big\}dz \right| \label{conta}\\ 
& = & \left|\intav{B_{r_n}} \int_0^1 \dfrac{d}{dt}\Big\{ f(tx + (1-t)y + z)g(tx + (1-t)y - z)\Big\} dt dz\right| \nonumber\\
& \leq & |y-x|\int_0^1\intav{B_{r_n}}\Big\{ |\nabla f(tx + (1-t)y + z)||g(tx + (1-t)y - z)| \nonumber\\
& &  \ \ \ \ \ \ \ \ \ \ \ \ \ \ \ \ \ + \ \ |f(tx + (1-t)y + z)||\nabla g(tx + (1-t)y - z)|\Big\}dz dt \nonumber\\
& \leq & |y-x| \int_0^1 \Big\{\mathcal{M}(|\nabla f|, g)(tx + (1-t)y) +  \mathcal{M}(f, |\nabla g|)(tx + (1-t)y)\Big\}dt \nonumber\\
& = & \int_{\overline{xy}}\Bigl\{\mathcal{M}(|\nabla f|, g) +  \mathcal{M}(f, |\nabla g|)\Bigr\} \, d\,\mathcal{H}^1, \nonumber
\end{eqnarray*}
we obtain 
\begin{equation}\label{abscont}
\left|\mathcal{M}(f,g)(x) - \mathcal{M}(f,g)(y)\right| \leq \int_{\overline{xy}} \Bigl\{\mathcal{M}(f,|\nabla g|) + \mathcal{M}(|\nabla f|,g)\Bigr\}\, d\,\mathcal{H}^1
\end{equation} 
for all $x,y \in \R^n$.

Now consider $f \in W^{1,p}(\R^n)$ and $g \in W^{1,q}(\R^n)$. Fix a vector $\nu \in S^{n-1}$ and consider sequences $\{f_j\}_{j=1}^{\infty}$ and $\{g_j\}_{j=1}^{\infty}$ of functions in $C^{\infty}_0$ such that
\begin{equation*}
  f_j \to f \ \ \textrm{in} \ \ W^{1,p}(\R^n) \ \ \ \ \textrm{and}  \ \ \ \ g_j \to g \ \ \textrm{in} \ \ W^{1,q}(\R^n).
\end{equation*}
From the continuity of the bilinear maximal operator in $L^p(\R^n) \times L^q(\R^n) \to L^r(\R^n)$ for $1/p + 1/q = 1/r$, $1<p,q<\infty$ and $r>1$ (here we can also include $n=1$ and $r=1$) we have
\begin{eqnarray}
\mathcal{M}(f_j, g_j) &\to& \mathcal{M}(f, g) \ \ \ \textrm{in} \ \ \ L^r(\R^n), \label{Sec2.1}\\
\mathcal{M}(f_j, |\nabla g_j|) &\to& \mathcal{M}(f, |\nabla g|) \ \ \ \textrm{in} \ \ \ L^r(\R^n), \label{Sec2.2}\\
\mathcal{M}(|\nabla f_j|, g_j) &\to& \mathcal{M}(|\nabla f|, g) \ \ \ \textrm{in} \ \ \ L^r(\R^n). \label{Sec2.3}
\end{eqnarray}
Using the fact that if $h_j \to h$ in $L^r(\R^n)$, then there is a subsequence such that for almost all lines $l$ parallel to $\nu$ the restriction of $h_j$ to $l$ converges in $L^r(l)$ to the restriction of $h$ to $l$, a standard approximation argument based on (\ref{abscont})-(\ref{Sec2.3}) gives
\begin{equation}\label{abscont2}
\left|\mathcal{M}(f,g)(x) - \mathcal{M}(f,g)(y)\right| \leq \int_{\overline{xy}} \Bigl\{\mathcal{M}(f,|\nabla g|) + \mathcal{M}(|\nabla f|,g)\Bigr\} \, d\,\mathcal{H}^1
\end{equation} 
almost everywhere on almost all lines parallel to $\nu$. This is sufficient to conclude that the weak derivative in the $\nu$-direction $D_{\nu}\mathcal{M}(f,g)(x)$ exists for almost every $x \in \R^n$ (cf. \cite[section 4.9]{EG}) and satisfies 
\begin{equation}
\left|D_{\nu}\mathcal{M}(f,g)(x)\right| \leq \mathcal{M}(f,|\nabla g|)(x) + \mathcal{M}(|\nabla f|,g)(x).
\end{equation}
Finally, taking the supremum over a countable and dense set of directions $\nu \in S^{n-1}$, we obtain 
\begin{equation}
|\nabla \mathcal{M}(f,g)(x)| \leq \mathcal{M}(f,|\nabla g|)(x) + \mathcal{M}(|\nabla f|,g)(x)
\end{equation}
for almost every $x \in \R^n$, which is (\ref{bound2}).

\section{Proof of Theorem 1: Continuity}%%%%%%%%%%%%%%%%%%%%%%%%%%%%%%%Section3%%%%%%%%%%%%%%%%%%%%%%%%%%%%%%%%%%%%%%%%%%%%%
Here we follow carefully the beautiful proof for the continuity in $W^{1,p}(\R^n)$ of the Hardy-Littlewood maximal operator in \cite{Lu}. Slight modifications are needed, and we simply adjust the notation to our context and quote the {\it lemmata} without proof in this sketch.

If $A \subset \R^n$ and $x \in \R^n$ we define
\begin{equation*}
 d(x,A) := \inf_{a\in A}|x-a| \ \ \ \textrm{and} \ \ \ A_{(\lambda)}:= \{x \in \R^n; d(x,A) \leq \lambda\} \ \ \textrm{for} \ \ \lambda \geq 0.
\end{equation*}
We recall that $W^{1,p}(\R^n)$ is endowed with the norm
\begin{equation*}
 \|f\|_{1,p} = \|f\|_{p} + \|\nabla f\|_{p}\ ,
\end{equation*}
where $\nabla f$ is the weak gradient of $f$. Denote $\|f\|_{p,A}$ for the $L^p$-norm of $\chi_{A} f$ for all measurable sets $A \subset \R^n$.

For $f \in L^p(\R^n)$ and $g \in L^q(\R^n)$ with $1<p,q<\infty$ and $1/p + 1/q = 1/r < 1$ ($1/r \leq 1$ if $n=1$), we define, for a fixed point $x \in \R^n$ the set of ``good'' radii $\mathcal{R}(f,g)(x)$ by
\begin{align}
\begin{split}
 \mathcal{R}(f,g)(x) = &\{r\geq 0; \, \mathcal{M}(f,g)(x) = \limsup_{r_k \to r} \intav{B_{r_k}} |f(x + y)g(x-y)|dy \\
 & \textrm{for some} \ \ r_k>0 \}.
\end{split}
\end{align}
From the definition it is clear that the set $\mathcal{R}(f,g)(x)$ is closed. If we define for each $x \in \R^n$ the function $u_x: [0,\infty) \to \R$ by
\begin{equation*}
 u_x(0) = |f(x)g(x)| \ \ \textrm{and} \ \ u_x(r) = \intav{B_r} |f(x + z)g(x-z)|dz \ \ \textrm{when} \ \ r \in (0,\infty),
\end{equation*}
then the functions $u_x$ are continuous on $(0,\infty)$ for all $x \in \R^n$. By an argument similar to the one that proves that almost every point is a Lebesgue point, we can also see that $u_x$ is continuous at $r=0$ for almost all $x \in \R^n$. By H\"older's inequality,
\begin{equation}
 u_x(r) \leq m(B_r)^{-\frac{1}{p} - \frac{1}{q}} \|f\|_p \|g\|_q\,  ,
\end{equation}
which proves that for almost every $x$ the function $u_x$ has at least one maximum point in $[0,\infty)$. Therefore the set $\mathcal{R}(f,g)(x)$ is nonempty and we have
\begin{eqnarray*}
&&\mathcal{M}(f,g)(x) =  \intav{B_r} |f(x + z)g(x-z)|dz \ \ \textrm{if} \ \ 0<  r\in \mathcal{R}(f,g)(x), \, \forall x\in \R^n, \\
&&\mathcal{M}(f,g)(x)  =  |f(x)g(x)| \ \ \textrm{for almost every} \ \ x \ \ \textrm{such that} \ \ 0 \in \mathcal{R}(f,g)(x).
\end{eqnarray*}
We refer now to \cite{Lu} for the ideas of the proofs for the following  {\it lemmata}.

\begin{lemma}[cf. Lemma 2.2 in \cite{Lu}]\label{lem2}

Suppose $f_j \to f$ in $L^p(\R^n)$, $g_j \to g$ in $L^q(\R^n)$ when $j \to \infty$. Then for all $R>0$ and $\lambda >0$ we have
\begin{equation*}
m(\{x \in B_R; \, \mathcal{R}(f_j,g_j)(x) \not\subset \mathcal{R}(f,g)(x)_{(\lambda)}\}) \to 0 \ \ \textrm{when} \ \ j \to \infty.
\end{equation*}

\end{lemma}
Let us introduce now some more notation. Let $e_i$ be one of the standard canonical vectors in $\R^n$. For $f \in L^p(\R^n)$ and $h > 0$, define
\begin{equation*}
 f^i_h(x) = \dfrac{f(x+he_i) - f(x)}{h} \ \ \ \textrm{and} \ \ \ f^i_{\tau(h)}(x) = f(x+he_i).
\end{equation*}
We know that for $p\geq 1$, $f^i_{\tau(h)} \to f$ in $L^p(\mathbb{R}^{n})$ when $h \to 0$, and if $f \in W^{1,p}(\R^n)$ we have $f^i_h \to D_if$ in $L^p(\R^n)$ when $h \to 0$, where $D_if$ denotes the partial derivative $\frac{\partial f}{\partial x_i}$.

The Hausdorff distance between two sets $A$ and $B$ is defined as
\begin{equation*}
\pi(A,B):= \inf \{\delta >0; \, A \subset B_{(\delta)} \ \ \textrm{and} \ \ B\subset A_{(\delta)}\}.
\end{equation*}
As a consequence of Lemma \ref{lem2} we have
\begin{lemma}[cf. Corollary 2.3 in \cite{Lu}]\label{lem3}
 Let $f \in L^p(\R^n)$, $g \in L^q(\R^n)$ with $1/p + 1/q = 1/r < 1$ $(1/r \leq 1$ if $n=1)$. Then for all $i$, $1\leq i \leq n$, $R>0$ and $\lambda>0$ we have
\begin{equation*}
 m(\{ x \in B_R; \, \pi(\mathcal{R}(f,g)(x), \mathcal{R}(f,g)(x +he_i)) > \lambda\}) \to 0 \ \ \textrm{when} \ \ h\to 0.
\end{equation*}

\end{lemma}

We are now in a position to state a formula for the derivative of the bilinear maximal function.

\begin{lemma}[cf. Theorem 3.1 in \cite{Lu}]\label{lem4}
Let $f \in W^{1,p}(\R^n)$, $g \in W^{1,q}(\R^n)$ with $1< p, q < \infty$ and $1/p + 1/q = 1/r < 1$ $(1/r \leq 1$ if $n=1)$. For almost all $x \in \R^n$ we have
\begin{eqnarray*}
D_i\mathcal{M}(f,g)(x) &=&  \intav{B_r} \Big\{D_i|f|(x + z)|g|(x-z) + |f|(x + z)D_i|g|(x-z)\Big\}dz \\
&&  \textrm{for all} \ \ 0<  r\in \mathcal{R}(f,g)(x),\\
\\
D_i\mathcal{M}(f,g)(x) & = & D_i|f|(x)|g|(x) + |f|(x)D_i|g|(x) \ \ \textrm{if} \ \ 0 \in \mathcal{R}(f,g)(x).
\end{eqnarray*}
\end{lemma}
With this machinery in hand, we are now able to prove the continuity of the bilinear maximal operator.

{\it Proof of the continuity}: Let $f_j \to f$ in $W^{1,p}(\R^n)$ and $g_j \to g$ in $W^{1,q}(\R^n)$. We must show that $\|\mathcal{M}(f_j,g_j) - \mathcal{M}(f,g)\|_{1,r} \to 0$. Since the bilinear maximal operator is sublinear, we know that $\|\mathcal{M}(f_j,g_j) - \mathcal{M}(f,g)\|_r \to 0$. This way, it suffices to prove that $\|D_i\mathcal{M}(f_j,g_j) - D_i\mathcal{M}(f,g)\|_r \to 0$ as $j \to \infty$ for all $1\leq i \leq n$. We may assume that the functions $f_j, g_j, f$ and $g$ are all nonnegative.

Fix $\epsilon >0$. Let us choose $R$  such that $\|2 \mathcal{M}(D_if,g)+2 \mathcal{M}(f,D_ig) \|_{r,C_1} < \epsilon$, where $C_1 = \R^n - B_R$. By absolute continuity, there exists $\eta >0$ such that $\|2 \mathcal{M}(D_if,g)+2 \mathcal{M}(f,D_ig)\|_{r,A} < \epsilon$ whenever $m(A) < \eta$ and $A$ is a measurable subset of $B_R$. Let us define
\begin{eqnarray*}
 v_x(f,g)(r) &=& \intav{B_r} \Big\{D_if(x + z)g(x-z) + f(x + z)D_ig(x-z)\Big\}dz, \\
v_x(f,g)(0) &=& D_if(x)g(x) + f(x)D_ig(x).
\end{eqnarray*}
As already observed, the functions $v_x(f,g)$ are continuous on $[0,\infty)$ for almost all $x \in \R^n$ and also converge to $0$ as $r \to \infty$. Therefore, for almost every $x$ the function $v_x(f,g)$ is uniformly continuous and we can find $\delta(x)>0$ such that 
\begin{equation*}
| v_x(f,g)(r_1) - v_x(f,g)(r_2) | < \dfrac{\epsilon}{m(B_R)^{\frac{1}{r}}} \ \ \textrm{whenever} \ \ |r_1 - r_2| < \delta(x).
\end{equation*}
We can write $B_R$ as
\begin{equation*}
 B_R = \left( \bigcup_{k=1}^{\infty}\Bigl\{x \in B_R;\, \delta(x) > \frac{1}{k}\Bigr\}\right) \cup \mathcal{N},
\end{equation*}
where $m(\mathcal{N}) =0$. From this we can choose $\delta >0$ such that 
\begin{align*}
\begin{split}
 m\Bigl(\Bigl\{x \in B_R; |v_x(f,g)(r_1) - v_x(f,g)(r_2)| & \geq \dfrac{\epsilon}{m(B_R)^{\frac{1}{r}}}\, , \\
& \textrm{for some} \, r_1, r_2 \, \textrm{with} \, |r_1 - r_2| <\delta \Bigr\}\Bigr) \\
:= m(C_2) < \frac{\eta}{2}.&
\end{split}
\end{align*}
Lemma \ref{lem2} says that we can find $j_0\in\N$ such that 
\begin{equation*}
 m(\{x \in B_R; \, \mathcal{R}(f_j,g_j)(x) \not\subset \mathcal{R}(f,g)(x)_{(\delta)}\}):= m(C^j) < \frac{\eta}{2} \ \ \textrm{when} \ \ j \geq j_0.
\end{equation*}
Fix $j \geq j_0$ and let $r_1 \in \mathcal{R}(f_j,g_j)(x)$, $r_2 \in \mathcal{R}(f,g)(x)$. From Lemma \ref{lem4} we have, for almost every $x \in \R^n$,
\begin{align}\label{est}
\begin{split}
|D_i&\mathcal{M}(f_j,g_j)(x) - D_i\mathcal{M}(f,g)(x)|  =  |v_x(f_j,g_j)(r_1) - v_x(f,g)(r_2)|\\
\\
&\leq |v_x(f_j,g_j)(r_1) - v_x(f_j,g)(r_1)| + |v_x(f_j,g)(r_1) - v_x(f,g)(r_1)|  \\
&  \ \ \ \ \ \ \ \ \ \ \ \ \ \ \ \ \ \ \ \ \ \ \ \  \ \ \ \ +  |v_x(f,g)(r_1) - v_x(f,g)(r_2)|  \\
\\
&\leq \mathcal{M}(D_if_j,g_j - g)(x) + \mathcal{M}(f_j, D_ig_j - D_ig)(x) \\
& \ \ \ \ \ \ \ \ \ \ \ \ \ \ \ \ \ \ + \mathcal{M}(D_if_j-D_if, g)(x) + \mathcal{M}(f_j - f, D_ig)(x) \\
& \ \ \ \ \ \ \ \ \ \ \ \ \ \ \ \ \ \ +  |v_x(f,g)(r_1) - v_x(f,g)(r_2)| .
\end{split}
\end{align}
If $x \notin C_1 \cup C_2 \cup C^j$ we can choose $r_1 \in \mathcal{R}(f_j,g_j)(x)$ and $r_2 \in \mathcal{R}(f,g)(x)$ such that $|r_1 - r_2| < \delta$. Our choice of $\delta$ then implies
\begin{equation}\label{est1}
 |v_x(f,g)(r_1) - v_x(f,g)(r_2)| < \dfrac{\epsilon}{m(B_R)^{\frac{1}{r}}}.
\end{equation}
If $x \in C_1 \cup C_2 \cup C^j$ we estimate
\begin{equation}\label{est2}
 |v_x(f,g)(r_1) - v_x(f,g)(r_2)| \leq 2\mathcal{M}(D_if,g)(x) + 2\mathcal{M}(f,D_ig)(x).
\end{equation}
Observe that $m(C_2 \cup C^j) <\eta$. Therefore, combining estimates (\ref{est1}) and (\ref{est2}) with the inequality (\ref{est}), we obtain
\begin{align*}
\begin{split}
\|D_i\mathcal{M}(f_j,g_j) -& D_i\mathcal{M}(f,g)\|_{r} \leq \| \mathcal{M}(D_if_j,g_j - g)\|_{r} + \|\mathcal{M}(f_j, D_ig_j - D_ig)\|_r \\
&+ \|\mathcal{M}(D_if_j-D_if, g)\|_r + \|\mathcal{M}(f_j - f, D_ig)\|_r \\
& +  \|\dfrac{\epsilon}{m(B_R)^{\frac{1}{r}}}\|_{r,B_R} + \|2\mathcal{M}(D_if,g) + 2\mathcal{M}(f,D_ig)\|_{r,C_1}\\
& + \|2\mathcal{M}(D_if,g) + 2\mathcal{M}(f,D_ig)\|_{r,C_2 \cup C^j}.
\end{split}
\end{align*}
The first four terms of the right-hand side converge to zero as $j \to \infty$ because of the boundedness of the bilinear maximal operator. The remaining terms are less than $\epsilon$ because of our choices for $R$ and $\eta$. Since $\epsilon>0$ was taken to be arbitrary, we conclude that $\|D_i\mathcal{M}(f_j,g_j) - D_i\mathcal{M}(f,g)\|_{r} \to 0$ as $j \to \infty$. This concludes the proof.

\section{Almost everywhere and weak convergence under the classical Hardy-Littlewood maximal operator}%%%%%%%%%%%%%%%%%%%%%%%%%%%%%%%%%%%%%%%%Section4%%%%%%%%%%%%%%%%%%%%%%%%%%%%%%%%%%%%%%%%%%%%%%%%%%%%%%%%%%%%%%%%%%%%%

We now turn our attention to the classical Hardy-Littlewood maximal operator to add some remarks to this theory. For $f:\R^n \to \R$ locally integrable it is defined as
\begin{equation}\label{Sec4.1}
 Mf(x) = \sup_{R>0}\dfrac{1}{m(B_R)}\int_{B_R(x)} |f(y)|dy,
\end{equation}
where $B_R(x)$ denotes the ball of radius $R$ centered in $x$. For $\Omega \subset \R^n$ a proper open subset of $\R^n$ and $f:\Omega \to \R$ we can define the local maximal operator at a point $x \in \Omega$ by
\begin{equation}\label{Sec4.2}
 M_{\Omega}f(x) = \sup_{0<R<\delta_x}\dfrac{1}{m(B_R)}\int_{B_R(x)} |f(y)|dy,
\end{equation}
where the supremum is taken over all radii $R$ such that $0 < R < \delta_x := \textrm{dist}(x,\partial \Omega)$. Both the global and the local maximal operators are known to be bounded from $L^p$ to $L^p$ when $p>1$. In this case the sublinearity of the operator implies continuity. As already mentioned in the Introduction of this paper, the operator $M$ defined in (\ref{Sec4.1}) is bounded and continuous from $W^{1,p}(\R^n)$ to $W^{1,p}(\R^n)$ (see \cite{K} and \cite{Lu}) and the local maximal operator $M_{\Omega}$ is bounded from $W^{1,p}(\Omega)$ to $W^{1,p}(\Omega)$ (see \cite{KL}).

We may ask ourselves if these classical maximal operators preserve other types of convergence, for instance pointwise convergence almost everywhere or weak convergence. The goal of this section is to settle the discussion about these issues, providing counterexamples and positive results in this direction.

\begin{example} 
The maximal operators $M: L^p(\R^n) \to L^p(\R^n)$ and $M_{\Omega}:L^p(\Omega) \to L^p(\Omega)$, for $p>1$, do not preserve pointwise convergence almost everywhere.
\end{example}
\begin{proof} This follows simply from the observation that 
\begin{equation}
M(f)(x) \geq C_{f} |x|^{-n} ~~\textrm{ whenever } \, |x|\geq 1,
\end{equation}
\noindent where $C_{f}:=\dfrac{||f||_{L^{1}(B_{1})}}{2^{n}\omega_{n}}$. Here $\omega_{n}$ is the volume of the unit ball $B_{1}\subset\R^{n}$. We consider the sequence
\begin{equation*}
u_k(x) = \dfrac{1}{m(B_{\frac{1}{k}})} \chi_{B_{\frac{1}{k}}}(x).
\end{equation*}
Clearly,  $u_k \to 0$ a.e. but $M(u_k) \not\rightarrow 0$ a.e. The argument for the local case is just a simple adaptation of this one.
\end{proof}

Issues about the stability of the weak convergence under nonlinear operators are much more interesting and have been studied in \cite{MT} for a certain class called Nemytskii nonlinearities with applications to differential equations in \cite{T}. Given an operator $T: E\to F$ between Banach spaces and $u_k \rightharpoonup u$ in $E$, the question is whether or not we have $T(u_k)\rightharpoonup T(u)$ in $F$ (in the affirmative case for all such sequences $\{u_{k}\}_{k\geq 1}$, we say that $T$ is sequentially weakly continuous). We show below a counterexample in this direction.

\begin{example}\label{exm6}
The maximal operators $M:L^p(\R^n) \to L^p(\R^n)$ and $M_{\Omega}: L^p(\Omega) \to L^p(\Omega)$, for $p >1$, are not sequentially weakly continuous.
\end{example}
\begin{proof}
We start with the local case. Let $\Omega = (-1, 1) \subset \R$ and consider the orthonormal system in $L^2(-1,1)$ given by $u_n(x) = \sin (2\pi nx)$, $n =1,2,3,...$. Therefore we have $u_n\rightharpoonup 0$ in $L^2(-1, 1)$, but we claim that $M_{\Omega}(u_n) \not\rightharpoonup 0$ in $L^2(-1,1)$.

To see this, let us fix a radius $r < \frac{1}{2}$ and consider the inner product
\begin{eqnarray*}
\langle 1, M_{\Omega}(u_n) \rangle_{L^2(-1, 1)} & = & \int_{-1}^{1} M_{\Omega}(u_n)(x)dx \geq \int_{-\frac{1}{2}}^{\frac{1}{1}} M_{\Omega}(u_n)(x)dx \\
& \geq & \int_{-\frac{1}{2}}^{\frac{1}{2}}\left( \intav{B_r} |u_n(x+y)|dy\right) dx\\
& =& \intav{B_r} \left(\int_{-\frac{1}{2}}^{\frac{1}{2}} |\sin(2\pi nx + 2\pi ny)|dx \right) dy\\
& = &  \intav{B_r} C dy = C ,
\end{eqnarray*}
 where $C >0$ is a constant. This proves our claim.

For the classical Hardy-Littlewood maximal operator we give the following counterexample in $L^2(\R)$:
\begin{equation*}
u_n(x) = \dfrac{\sin(2\pi nx)}{1 + x^2}.
\end{equation*}
We have $u_n \rightharpoonup 0$ in $L^{2}(\R)$ as a consequence of the Riemann-Lebesgue Lemma, but we claim that $M(u_n) \not\rightharpoonup 0$ in $L^{2}(\R)$. To verify this, fix $1>r>0$ and observe that
\begin{align}\label{Sec4.4}
\begin{split}
 \langle M(u_n), \frac{1}{1+x^2} \rangle _{L^2(\R)} &= \int_{\R} M(u_n)(x) \frac{1}{1+x^2}dx \\
& \geq  \int_{\R} \left( \intav{B_r} |u_n(x+y)|dy \right)\frac{1}{1+x^2}dx  \\
& =  \intav{B_r} \left(\int_{\R} \frac{|\sin(2\pi n (x+y)|}{(1 + (x+y)^2)(1+x^2)}dx \right) dy \\
& \geq  \intav{B_r} C dy = C,
\end{split}
\end{align}
where $C>0$ is a constant. This finishes the proof of our claim.
\end{proof}
We remark that the previous example also proves that $M:L^{\infty}(\R^n) \to L^{\infty}(\R^n)$ and $M_{\Omega}: L^{\infty}(\Omega) \to L^{\infty}(\Omega)$ {\it are not sequentially weakly star continuous}.

It is interesting to compare the previous example with the following positive results in Sobolev spaces. 

\begin{proposition}\label{prop8}
Suppose $\Omega$ is a bounded domain with Lipschitz boundary. Then, the local maximal operator $M_{\Omega}:W^{1,p}(\Omega) \to W^{1,p}(\Omega)$ is sequentially weakly continuous for $p>1$.
\end{proposition}
\begin{proof}
Let $f_j \rightharpoonup f$ in $W^{1,p}(\Omega)$. Since $M_{\Omega}:W^{1,p}(\Omega) \to W^{1,p}(\Omega)$ is a bounded operator, the sequence $M_{\Omega}(f_j)$ must admit a weakly convergent subsequence, by reflexivity. This way, we can assume
\begin{equation*}
M_{\Omega}(f_j) \rightharpoonup g \ \ \textrm{in} \ \  W^{1,p}(\Omega).
\end{equation*}
 By the compactness of the Sobolev embedding $W^{1,p}(\Omega)\hookrightarrow L^{p}(\Omega)$ and the continuity in $L^{p}(\Omega)$ of the local maximal operator, we have 
\begin{equation*}
M_{\Omega}(f_j) \rightarrow  M_{\Omega}(f)\ \ \textrm{in} \ \  L^p(\Omega).
\end{equation*}
In particular, $M_{\Omega}(f)=g$ and this finishes the proof. 
\end{proof}

\begin{theorem}\label{thm9}  
Let $1<p<\infty$ and suppose $u_{k}\rightharpoonup  u$ in $W^{1,p}(\R^{n})$. There exists a subsequence $M(u_{k_{j}})\to M(u)~~ a.e.$ in $\R^{n}$.
\end{theorem}

\begin{proof} By the sublinearity of the maximal operator, it is enough to prove the case where $u\equiv 0$. Let us consider $B=B_{L}(0)$, where $L>0$.

First, we observe that if $f\in L^{p}(\R^{n})$, there exists a universal $C>0$ (depending only on the dimension $n$) such that 
\begin{equation*}
\intav{B_{R}(x)} |f(y)|dy \leq CR^{-\frac{n}{p}}||f||_{L^{p}(\R^n)}.
\end{equation*}
For each $m = 1,2,3,...$ we can take $R_{m} >0$ large enough so that, for every $x\in\mathbb{R}^{n}$,

\begin{equation}\label{pequeno}
\intav{B_{R}(x)} |u_{k}(y)|dy \leq \frac{1}{m} \hspace{.2cm} \textrm{ for all } k\in\N \hspace{.2cm} \textrm{ whenever } R\geq R_{m}.
\end{equation}
Let us consider now $B^{\star}_m:=B_{L+2R_{m}}(0)$ and the local maximal operator with respect to this ball.  For any $x\in B$, if $\delta_{x}=\textrm{dist}(x,\partial B^{\star}_m)$, we have by the estimate (\ref{pequeno})

\begin{align}\label{splitting}
\begin{split}
  M(u_{k})(x)& =\max\Big\{\sup\limits_{0<R<\delta_{x}} \intav{B_{R}(x)} |u_{k}(y)|dy, \sup\limits_{R>\delta_{x}} \intav{B_{R}(x)} |u_{k}(y)|dy \Big\} \\
& \leq \max\Big\{M_{B^{\star}_m}(u_{k})(x), \frac{1}{m} \Big\} .
 \end{split}
 \end{align}
Since $W^{1,p}(\R^{n}) \hookrightarrow W^{1,p}(B^{\star}_m)$ continuously and 
$W^{1,p}(B^{\star}_m) \hookrightarrow L^{p}(B^{\star}_m)$ compactly, by the continuity of the local maximal operator $M_{B^{\star}_m}$, $M_{B^{\star}_m}(u_{k}) \to 0$ in $ L^{p}(B^{\star}_m)$. Therefore, there is a subsequence $ M_{B^{\star}_m}(u^m_{k_{j}})\to 0~~$ a.e. in $B^{\star}_m$. From (\ref{splitting}) we conclude that 
$$ \limsup_{j\to \infty} M(u^m_{k_{j}})(x) \leq \frac{1}{m} \ \ \textrm{ a.e. in } B. $$ 
Using the Cantor diagonal argument we can find a subsequence $\{u_{k_j}\}$ such that 
$$M(u_{k_{j}})(x) \rightarrow 0 \ \ \textrm{ a.e. in } B. $$ 
Since the original ball $B$ was arbitrary, we can use once more the Cantor diagonal argument applied to $\mathbb{R}^{n}=\bigcup\limits_{n=0}^{\infty} B_{n}(0)$ to conclude the theorem. 
\end{proof}

\begin{corollary} Assume $1<p<\infty$. The maximal operator $M:W^{1,p}(\R^{n})\to W^{1,p}(\mathbb{R}^{n})$ is sequentially weakly continuous. 

\end{corollary}

\begin{proof} The proof is similar to the proof for the local maximal operator given in Proposition \ref{prop8}, with the help of the previous theorem. Let $u_{k}\rightharpoonup u$ in $W^{1,p}(\mathbb{R}^{n})$. By the boundedness of the maximal operator in $W^{1,p}(\R^n)$, we can assume $M(u_{k})\rightharpoonup g$ in $W^{1,p}(\mathbb{R}^{n})$. By the previous theorem, there exists a subsequence $ M(u_{k_{j}})\to M(u)$ a.e. in $\R^n$. This is sufficient to conclude that $M(u) = g$.
\end{proof}

We observe that Theorem \ref{thm9} is optimal in the sense that one cannot replace the weak convergence in $W^{1,p}(\R^n)$ for weak convergence in $L^p(\R^n)$. Example \ref{exm6} above presents a sequence $u_k \rightharpoonup 0$ in $L^2(\R)$ such that $M(u_k) \not\rightarrow 0$ a.e.

Finally, we point out that Theorem \ref{thm9} and its corollary are also optimal in the right-hand side. We present an example showing that the maximal operator is not compact in the sense that it does not map weakly convergent sequences into strongly convergent sequences.

\begin{example} 
For $1<p<\infty$, the maximal operators $M:W^{1,p}(\R^n) \to L^p(\R^n)$ and $M_{\Omega}:W^{1,p}(\Omega) \to L^p(\Omega)$ are not compact.
\end{example}
\begin{proof}
For the local case, consider the sequence of disjoint balls $B_{k}:=B_{1/2}(ke_{1})$, where $e_{1}=(1,0,\dots,0)\in\R^{n}$, put $\Omega = \bigcup_{k=1}^{\infty} B_k$ and take the sequence of functions $u_{k}:= m(B_{k})^{-1/p}\chi_{B_{k}}\in C^{\infty}(\Omega)$. For the global case let $u \in C^{\infty}_0(\R^n)$ and consider the sequence $u_k(x) = u(x-k)$.
\end{proof}

\section*{Acknowledgments}
We would like to thank William Beckner for the motivation and valuable suggestions. We are also grateful to Michael Lacey for the remark after Theorem 1.

\bibliographystyle{amsplain}

\begin{thebibliography}{99}
\bibitem{AP} J.M. Aldaz and J. Pérez Lázaro, {\it Functions of bounded variation, the derivative of the one dimensional maximal function, and applications to inequalities},  Trans. Amer. Math. Soc. {\bf 359} (2007), no. 5, 2443--2461.
\bibitem{EG} L.C. Evans and R.F. Gariepy, {\it Measure theory and fine properties of functions}, Stud. Adv. Math., CRC, Boca Raton, FL, 1992.
\bibitem{HO} P. Haj\l asz and J. Onninen, {\it On boundedness of maximal functions in Sobolev spaces},  Ann. Acad. Sci. Fenn. Math.  {\bf 29}  (2004),  no. 1, 167--176.
\bibitem{K} J. Kinnunen, {\it The Hardy-Littlewood maximal operator of a Sobolev function}, Israel J. Math. {\bf 100} (1997), 117--124.
\bibitem{KL} J. Kinnunen and P. Lindqvist, {\it The derivative of the maximal function},  J. Reine Angew. Math.  {\bf 503}  (1998), 161--167.
\bibitem{KS} J. Kinnunen and E. Saksman, {\it Regularity of the fractional maximal function}, Bull. London Math. Soc. {\bf 35} (2003), no. 4, 529--535.
\bibitem{L} M. Lacey, {\it The bilinear maximal function map into $L^p$ for $2/3 < p \leq 1$}, Ann. of Math. {\bf 151} (2000), 35--57.
\bibitem{Lu} H. Luiro, {\it Continuity of the maximal operator in Sobolev spaces},  Proc. Amer. Math. Soc.  {\bf 135}  (2007),  no. 1, 243--251.
\bibitem{MT} D. Moreira and E. Teixeira, {\it On the behavior of weak convergence under nonlinearities and applications}, Proc. Amer. Math. Soc. {\bf 133} (2005), no. 6, 1647--1656.
\bibitem{T} E. Teixeira, {\it Strong solutions for differential equations in abstract spaces}, J. Differential Equations {\bf 214} (2005), no. 1, 65--91.
\end{thebibliography}

\end{document}